\documentclass[reqno]{conm-p-l}
\usepackage{amssymb,accents}
\input{diagrams}
\def\noi{\noindent}

\newtheorem{Thm}{Theorem}[section]
\newtheorem{Lm}[Thm]{Lemma}
\newtheorem{Prop}[Thm]{Proposition}
\newtheorem{Cor}[Thm]{Corollary}
\newtheorem{Def}[Thm]{Definition}
\newtheorem{state}{Definition}
\theoremstyle{remark}
\newtheorem{Rem}[Thm]{Remark}

\numberwithin{equation}{section}

\def\cal{\mathcal}
\def\Bbb{\mathbb}
\def\mf{\mathfrak}

\def\<{\langle}
\def\>{\rangle}

\def\a{\alpha}
\def\b{\beta}
\def\d{\delta}
\def\D{\Delta}
\def\th{\theta}

\def\l{\lambda}
\def\L{\Lambda}

\def\Re{\Bbb R}
\def\F{\Bbb F}

\def\Z{\Bbb Z}

\def\T{\cal T}

\def\H{\cal H}

\def\A{\cal A}
\def\q{\cal C}
\def\D{\cal D}

\def\G{\mf h}

\def\W{\ring W}
\def\w{\ring w}
\def\RR{\ring R}
\def\Q{\ring Q}

\def\Gc{\ring {\mf h}}

\begin{document}
\title{Triple groups and Cherednik algebras}
\author{Bogdan Ion}
\address{Department of Mathematics, University of Michigan, Ann
Arbor, MI 48109}
\email{bogdion@umich.edu}
\author{Siddhartha Sahi}
\address{Department of Mathematics,
Rutgers University, New Brunswick, NJ 08903 }
\email{sahi@math.rutgers.edu}
\subjclass[2000]{Primary 20C08; Secondary 33D52}
\copyrightinfo{2005}{American Mathematical Society}

\begin{abstract}
The goal of this paper is to define a new class of objects which
we call triple
 groups and to relate them with
Cherednik's double affine Hecke algebras. This has as immediate
consequences new descriptions of double affine Weyl and Artin
groups, the double affine Hecke algebras as well as the
corresponding elliptic objects. From the new descriptions we
recover results of Cherednik on automorphisms of double affine
Hecke algebras.
\end{abstract}
\maketitle

%%%%%%%%%%%%%%%%%%%%%%%%%%%%%%%%%%%%%%%%%%%%%%%%%%%%%%%%%%%%%%%%%%%%%%%%%%%%%
%%%%%%%%%%%%%%%%%%%%%%%%%%%%%%%%%%%%%%%%%%%%%%%%%%%%%%%%%%%%%%%%%%%%%%%%%%%%%
\thispagestyle{empty}
\section*{Introduction}

The classical theory of root systems associated with finite
dimensional semisimple Lie algebras admits two generalizations.
The first one, the Kac--Moody theory, originated in the effort to
extend the classical theory to include certain classes of infinite
dimensional Lie algebras as, for example, the simple infinite
dimensional Lie algebras of vector fields on finite dimensional
spaces as classified by  Cartan at the beginning of the 20th
century. The theory, which was initiated in the mid 1960's by V.
Kac and R. Moody, grew to be extremely rich and with abundant
connections with diverse areas of mathematics. The second
generalization, axiomatized by K. Saito, emerged from the theory
of semi-universal deformations of simply elliptic singularities
and the structure of Milnor lattices attached to these
singularities. This second theory is still in the very early
stages of development from a representation--theoretical point of
view, at present extending only the classical theory of root
systems. These root systems were called by Saito $2$--extended
root systems or elliptic root systems. Saito defined in fact more
general classes of root systems, $k$--extended root systems ($k$
being any positive integer), but for values of $k$ higher than
three these do not have any geometrical or
representation--theoretical relevance yet. The problem of
attaching corresponding notions of Lie algebras and groups to
$k$--extended root systems is wide open although in the elliptic
case some constructions exist in the literature.

To briefly compare the two theories note that they include the
classical theory and they agree on one other class of root
systems: the 1-extended root systems coincide with affine root
systems in Kac--Moody theory. With this exception, they produce
quite different outcomes. To mention a few differences note that
the Weyl groups associated to Kac--Moody root systems are always
Coxeter groups, but the corresponding groups associated to
$k$--extended root systems are not. Also, the Weyl group of a
$k$--extended root system does not act properly discontinuously
anywhere in its reflection representation and as a consequence we
do not have a Tits cone, and no natural notion of positive root
system or Dynkin diagram. Saito was nevertheless able to classify
the irreducible elliptic root systems parameterizing them by what
he called ``elliptic Dynkin diagrams''.

Since, as mentioned above, the Weyl groups  associated to elliptic
root systems (also called elliptic Weyl groups) are not Coxeter
groups, a presentation of them in terms of elliptic Dynkin
diagrams was one of the first questions to be asked. A candidate
for such a presentation was exposed in \cite{saito2}. The
presentation is complicated and does not bring any light into the
structure of the elliptic Weyl group. Another problem along these
lines which was asked in \cite{saito2} is the description of the
fundamental group of the complement of the discriminant of the
semi--universal deformation of simply--elliptic singularities: the
elliptic Artin group. In fact, the elliptic Weyl groups appeared
before of the formal definition of elliptic root systems in the
work of Looijenga \cite{looijenga} who considered  extended
Coxeter groups. Under this terminology the elliptic Weyl groups
would be called extended affine Weyl groups. Also, the elliptic
Artin groups were described (in the non-exceptional cases) by
generators and relations by H. van der Lek \cite{lek} under the
name of extended Artin groups. The elliptic Artin group was
recognized by Cherednik \cite{c0} as being closely related to his
double affine Artin groups and Hecke algebras.

At this point a broad picture emerges since Cherednik's work has
deep connections with quantum many--body problems, conformal field
theory and  also with representation theory via Macdonald theory,
which encloses information about real and $p$--adic groups and
also about affine Kac--Moody groups. There are also connections
with harmonic analysis, $q$--Riemann zeta functions,
Gauss--Selberg sums and more. The area is in full development and
interesting connections are still to be made. In all the
applications the structure of double affine Artin groups and Hecke
algebras is especially important.

The first results on the structure of double affine Hecke algebras
appeared in Cherednik's work on Macdonald's conjectures. For the
proof of the constant term conjecture only basic properties of
Hecke algebras are used, but for the proof of the
evaluation--duality conjecture a deeper result was needed: the
existence of a special involution on double affine Hecke algebras
called the duality involution. The existence of this
involution is certainly non--obvious, see for example \cite{ion},
\cite[Chapter 3]{macbook} for a proof. Cherednik discovered that
in fact this is part of an action of ${\rm GL}(2,\Z)$ as outer
automorphisms of double affine Artin groups. This allows, for
example, representations of ${\rm PSL}(2,\Z)$ for which the matrix
coefficients are expressed in terms of special values of Macdonald
polynomials at roots of unity. Let us mention that in the case of
a root system of type $A_n$ such actions appeared before in the
work of A. Kirillov, Jr \cite{kirillov} on modular tensor
categories and quantum groups at roots of unity. In this case the
two approaches are most probably connected via the equivalence
(see \cite{kl}) between certain categories of representations of
quantum groups and affine Kac--Moody algebras. For the root system
of type $A_n$ the picture seems to be even richer since ${\rm
PSL}(2,\Z)$ (which is the mapping class group of the torus with
one marked point) can be replaced with the mapping class group of
any two dimensional surface with marked points.

The main goal of this paper is to give a completely new
description of double affine Hecke algebras, Artin groups and
double affine Weyl groups, and consequently of the corresponding
elliptic objects. Our description is much simpler than the
existing ones and has the virtue of making the existence of the
difference Fourier transform and of the ${\rm PSL}(2,\Z)$ action
mentioned above simple consequences. All our proofs are at the
level of Artin groups, since any other result is then obtained by
passing to various quotients. It turns out that in fact all these
symmetries of double affine Artin groups are descending from those
of a new object, which we call triple %%%%%%%%%%%%%%%%%%%%%%%%%%%%%%%%%%%%%%
% affine Artin
 group, of which the double affine Artin group is a quotient. A brief
description of our results follows. The reader should refer to
Section \ref{prelims} for notation and conventions.
\begin{state}\label{maindef}
Let $A$ be an irreducible affine Cartan matrix subject to our
restriction and $S(A)$ its Dynkin diagram. The triple %%%%%%%%%%%%%%%%%%%%%%%
%affine Artin
group $\A$ is given by generators and relations as follows:

\underline{Generators}: one generator $\T_i$ for each node, with
the exception of the affine node for which we have three
generators $\T_{01}$, $\T_{02}$ and $\T_{03}$. We use the notation
$\D=\T_{01}\T_{02}\T_{03}$.

\underline{Relations}: a)
 Braid relations for each pair of generators
associated to any pair of distinct nodes (note that there are
three generators associated to the affine node).

\hspace{1.52cm} b) If the affine node is connected with the node
$\a$ by a single lace, the elements $\T_{01}$, $\T_{02}$,
$\T_{03}$, $\T_{01}\T_{02}\T_{01}^{-1}$ and
$\T_{03}^{-1}\T_{02}\T_{03}$ satisfy the single lace Coxeter
relation with $\D^k \T_\a \D^{-k}$ for all integers $k$.

\hspace{1.52cm} c) If there are double laces connecting the affine
node with the node $\a$ the following relation holds for any
$1\leq i<j\leq3$
\begin{equation}\label{ellbraid}
\T_{0i}\T_\a^{-1}\T_{0j}\T_\a=\T_\a^{-1}\T_{0j}\T_\a\T_{0i}
\end{equation}
\end{state}

Of course we can define the corresponding triple %%%%%%%%%%%%%%%%%%%%%%%%%%
 Weyl group $\cal
W$ by further asking that the generators have order two. All the
facts we discuss will have as an immediate consequence
corresponding results for the triple %%%%%%%%%%%%%%%%%%%%%%%%%%%%%%%%%%%%%%
Weyl group.

A few comments on this definition may be useful at this stage. We
would like to stress that the above definition is not symmetric in
the generators $\T_{0i}$ which means that the symmetric group on
three letters will not act on the triple
 group.
However, we will prove in Theorem \ref{tripleartinauto} and
Theorem \ref{triplefaith} that the braid group on three letters
acts faithfully as automorphisms of the triple %%%%%%%%%%%%%%%%%%%%%%%%%%%%%
 group. What is
completely obvious from this Definition is the existence of an
{\sl anti-involution} which fixes all the generators except
$\T_{01}$ and $\T_{03}$ which are interchanged. Also, we like to
note that in Definition \ref{maindef}, b) we ask the relations for
the element $\T_{01}\T_{02}\T_{01}^{-1}$ merely for symmetry. They
follow from the ones for $\T_{03}^{-1}\T_{02}\T_{03}$ by observing
that
$$
\T_{01}\T_{02}\T_{01}^{-1}=\D \T_{03}^{-1}\T_{02}\T_{03}\D^{-1}.
$$

Our main result, Theorem \ref{main}, shows that the double affine
Artin group is a quotient of the triple
 group. In
fact, in this quotient  the relations in Definition \ref{maindef},
b) will become redundant, allowing us to give a new description of
double affine Artin groups and Hecke algebras.

For double affine Artin groups the descent of the above
anti--involution, whose existence is a trivial consequence of our
new description, gives (by composition with the anti-involution
which sends all the generators to their inverses) the duality
involution responsible for the difference Fourier transform. Note
that previous results explaining in an uniform way the existence
of this involution required topological arguments \cite{c0},
\cite{ion}. There is also an algebraic proof of this result \cite{macbook}. This proof
proceeds by
carefully  checking  of all the necessary relations between the
generators in the Cherednik presentation and
requires special considerations for some types of root systems.

Also, the above action  of the braid group on three letters
descends to double affine Artin groups and Hecke algebras. Here we
rediscover results of Cherednik which explain that there is a
morphism from the modular group to the group of outer
automorphisms of an double affine Artin group. These results are
described in Theorem \ref{involution}, Theorem
\ref{doubleartinauto} and Corollary \ref{cheredthm}. As another
consequence we obtain descriptions of elliptic Weyl groups and
elliptic Artin groups for certain types of elliptic root systems.

There are interesting quotients of  the triple %%%%%%%%%%%%%%%%%%%%%%%%%%%%%%
groups of which the double affine Artin groups are quotients but
which, unlike the triple %%%%%%%%%%%%%%%%%%%%%%%%%%%%%%%%%%%%%%%%%%%%%%%%%%%%
groups, are described by a finite number of relations. These
quotients also inherit the anti--involution and the action of the
braid group on three letters. It is not yet clear if any of these
quotients (or the triple %%%%%%%%%%%%%%%%%%%%%%%%%%%%%%%%%%%%%%%%%%%%%%%%%%%%
group itself) have a topological interpretation.

%%%%%%%%%%%%%%%%%%%%%%%%%%%%%%%%%%%%%%%%%%%%%%%%%%%%%%%%%%%%%%%%%%%%%%%%%%%%%
%%%%%%%%%%%%%%%%%%%%%%%%%%%%%%%%%%%%%%%%%%%%%%%%%%%%%%%%%%%%%%%%%%%%%%%%%%%%%
%%%%%%%%%%%%%%%%%%%%%%%%%%%%%%%%%%%%%%%%%%%%%%%%%%%%%%%%%%%%%%%%%%%%%%%%%%%%%
\section{Preliminaries}\label{prelims}

\subsection{Notation and conventions}

For the most part we adhere to the notation in \cite{kac}. Let
$A=(a_{jk})_{0\leq j,k\leq n}$ be an irreducible \emph{affine}
Cartan matrix, $S(A)$ the Dynkin diagram and  $(a_0,\dots, a_n)$
the numerical labels of $S(A)$ in Table Aff from
\cite[p.48-49]{kac}. Note that we consider that the nodes $i$ and $j$ are connected
in the Dynkin diagram by $a_{ij}a_{ji}$ laces. Unless  $A=A^{(1)}_1$ this produces the
same diagrams as in \cite{kac} (We thank Craig Johnson for pointing this out).
We denote by $(a_0^\vee,\dots, a_n^\vee)$ the
labels of the Dynkin diagram $S(A^t)$ of the dual algebra which is
obtained from $S(A)$ by reversing the direction of all arrows and
keeping the same enumeration of the vertices.

Let $({\G}, R, R^{\vee})$ be a realization of the Cartan matrix
$A$ and let $(\Gc, \RR, \RR^{\vee})$ be the associated finite root
system (which is a realization of the Cartan matrix $\ring A =
(a_{jk})_{1\leq j,k\leq n}$). If we denote by $\{\a_j\}_{0\leq
j\leq n}$ a basis of $R$ such that $\{\a_j\}_{1\leq j\leq n}$ is a
basis of $\RR$ we have the following description
$$
{\G}^*={\Gc^*} + {\mathbb R}\delta + {\mathbb R}{\Lambda}_0\ ,
$$
where $\d=\sum_{j=0}^n a_j\a_j$. The vector space ${\G}^*$ has a
canonical scalar product defined as follows
$$
(\a_j,\a_k):=d_j^{-1}a_{jk}\ ,\ \ \ \
(\L_0,\a_j):=\d_{j,0}a_0^{-1}\ \ \ \text{and}\ \ \ (\L_0,\L_0):=0,
$$
with $d_j:= a_ja_j^{{\vee}-1}$ and $\d_{j,0}$  Kronecker's delta.
The simple coroots are $\{\a_j^\vee:=d_j\a_j\}_{0\leq j\leq n}$.
The lattice $\Q =\oplus_{j=1}^n\Z\a_j$ is the root lattice of
$\RR$ and if the affine Cartan matrix is not $A_{2n}^{(2)}$ the
lattice $Q=\oplus_{j=0}^n\Z\a_j=\Q\oplus \Z\d$ is the root lattice
of $R$. For the affine Cartan matrix $A_{2n}^{(2)}$ we will have
to replace the root lattice $\Q$ by the weight lattice $P$ and to
keep in mind that the root lattice of $R$ is
$Q=\oplus_{j=0}^n\Z\a_j=P\oplus \frac{1}{2}\Z\d$.

\bigskip
\begin{quotation}\sl
To keep the notation as simple as possible we agree to use the
same symbol $\Q$ to refer to the root lattice or the weight
lattice depending on the affine root system as explained above.
\end{quotation}
\bigskip

 Although the ideas we will present here apply in general, the objects
associated to the irreducible affine Cartan matrices $B_n^{(1)}$,
$C_n^{(1)}$, $F_4^{(1)}$ and $G_2^{(1)}$ require special treatment
and we will not study them here. Their special behavior is a
consequence of the fact that they lack certain types of symmetries
which are otherwise abundant.
\bigskip
\begin{quotation}
{\sl
 For the rest of the paper we assume our irreducible
affine root system to be such that the affine simple root $\a_0$
is short (this  includes of course the case when all the roots
have the same length). The exceptions to this condition are:
$B_n^{(1)}$, $C_n^{(1)}$, $F_4^{(1)}$ and $G_2^{(1)}$. }
\end{quotation}
\bigskip

For any root system denote by $\a$ the simple root corresponding
to the node in the Dynkin diagram which is connected to node
associated to the affine simple root $\a_0$. If there are two such
nodes (which is the case for $S(A_n^{(1)})$) we choose one of
them. By $l_0$ we denote the number of laces by which the nodes
corresponding to $\a_0$ and $\a$ are connected. Note that $l_0$ is
always 1, except for $A= B_n^{(2)}, A_{2n}^{(2)}$  for which it
takes the value 2.
%%%%%%%%%%%%%%%%%%%%%%%%%%%%%%%%%%%%%%%%%%%%%%%%%%%%%%%%%%%%%%%%%%%%%%%%%%%%%
%%%%%%%%%%%%%%%%%%%%%%%%%%%%%%%%%%%%%%%%%%%%%%%%%%%%%%%%%%%%%%%%%%%%%%%%%%%%%

\subsection{Affine and double affine Weyl groups}
Given $\a\in R$, $x\in \G^*$ let
$$
s_\a(x):=x-\frac{2(x,\a)}{(\a,\a)}\a\ .
$$
The {\sl affine Weyl group} $W$ is the subgroup of ${\rm
GL}(\G^*)$ generated by all $s_\a$ (the simple reflections
$s_j=s_{\a_j}$ are enough). The {\sl finite Weyl group} $\W$ is
the subgroup generated by $s_1,\dots,s_n$. Both the finite and the
affine Weyl group are Coxeter groups and they can be abstractly
defined as generated by $s_1,\dots,s_n$, respectively
$s_0,\dots,s_n$, and some relations. These relations are called
Coxeter relations and they are of two types:
\begin{enumerate}
\item[a)] reflection relations: $s_j^2=1$; \item[b)] braid
relations: $s_is_j\cdots =s_js_i\cdots $ (there are $m_{ij}$
factors on each side, $m_{ij}$ being equal to $2,3,4,6$ if the
number of laces connecting the corresponding nodes in the Dynkin
diagram is $0,1,2,3$ respectively).
\end{enumerate}

\medskip
The {\sl double affine Weyl group} $\tilde W$ is defined to be the
semidirect product $W\ltimes Q$ of the affine Weyl group and the
lattice $Q$ (regarded as an abelian group with elements $\tau_\b$,
where $\b$ is a root), the affine Weyl group acting on the root
lattice as follows
$$
w\tau_\b w^{-1}=\tau_{w(\b)}.
$$
This group is the hyperbolic extension of an elliptic Weyl group
which is by definition the Weyl group associated with an elliptic
root system (see \cite{saito2} for definitions). The cases we
consider cover all  elliptic root systems with equal labels
$(p,p)$ (some of these will have isomorphic elliptic Weyl groups).
It also has a presentation with generators and relations (called
elliptic Coxeter relations). We refer the reader to \cite{saito2}
for the details. We only note here that the above mentioned
presentation is considerably more complex than the one for Coxeter
groups: for the infinite series it involves at least $2n-2$
generators and the genuine elliptic Coxeter relations (i.e. which
are not Coxeter relations) are relations among groups of three or
four of the generators.

The affine Weyl group $W$ can also be presented as a semidirect
product in the following way: it is the semidirect product of $\W$
and the lattice $\Q$ (regarded as an abelian group with elements
$\l_\mu$, where $\mu$ is in $\Q$), the finite Weyl group acting on
the root lattice as follows
$$
\w\l_\mu \w^{-1}=\l_{\w(\mu)}.
$$
Using this description we immediately see that the double affine
Weyl group could be described as follows.
\begin{Prop}\label{wdef}
The double affine Weyl group is the group generated by the finite
Weyl group $\W$, two lattices $\{\l_\mu\}_{\mu\in \Q}$,
$\{\tau_\b\}_{\b\in \Q}$ and an element $\tau_{a_0^{-1}\d}$ with
the following relations:
\begin{enumerate}
\item[(i)] $\w\l_\mu \w^{-1}=\l_{\w(\mu)}$ and $\w\tau_\b
\w^{-1}=\tau_{\w(\b)}$ for any $\w$ in the finite Weyl group and
any $\mu, \b$ in the root lattice; \item[(ii)]
$\l_\mu\tau_\b=\tau_\b\l_\mu\tau_{-(\b,\mu)\d}$; \item[(iii)]
$\tau_{a_0^{-1}\d}$ is central.
\end{enumerate}
\end{Prop}
To elucidate the alluded relation between the affine Weyl group
and the elliptic Weyl group we give the following definition.
\begin{Def}
The elliptic Weyl group is the factor group of the double affine
Weyl group by the group generated by $\tau_{a_0^{-1}\d}$.
\end{Def}

For $r$ a real number, $\G^*_r=\{ x\in\G\ ;\ (x,\d)=r\}$ is the
level $r$ of $\G^*$. We have
$$
\G^*_r=\G^*_0+r\L_0=\Gc^*+{\Bbb R}\d+r\L_0\ .
$$
The action of $W$ preserves each $\G^*_r$ and  we can identify
each $\G^*_r$ canonically with $\G^*_0$ and obtain an (affine)
action of $W$ on $\G^*_0$. For example, the level zero action of
$s_0$ and $\l_\mu$ on $\G^*_0$ is
\begin{eqnarray*}
s_0(x)     & = & s_\th(x)+(x,\th)a_0^{-1}\d\ ,\\
\l_\mu(x)   & = & x - (x,\mu)\d\ ,
\end{eqnarray*}
and the full action of the same elements on $\G^*$ is
\begin{eqnarray*}
s_0\<x\>   & = & s_\th(x)+(x,\th)a_0^{-1}\d-(x,\d)\a_0\ ,\\
\l_\mu\<x\> & = & x  - (x,\mu)\d+ (x,\d)(\mu
-\frac{1}{2}|\mu|^2\d) \ ,
\end{eqnarray*}
where we denoted by $\th=\delta -a_0\a_0$.
%%%%%%%%%%%%%%%%%%%%%%%%%%%%%%%%%%%%%%%%%%%%%%%%%%%%%%%%%%%%%%%%%%%%%%%%%%%%%
%%%%%%%%%%%%%%%%%%%%%%%%%%%%%%%%%%%%%%%%%%%%%%%%%%%%%%%%%%%%%%%%%%%%%%%%%%%%%
%%%%%%%%%%%%%%%%%%%%%%%%%%%%%%%%%%%%%%%%%%%%%%%%%%%%%%%%%%%%%%%%%%%%%%%%%%%%%

\subsection{Artin groups and Hecke algebras}
\bigskip
To any Coxeter group we can associate its Artin group, the group
defined with the same generators which satisfy only the  {\sl
braid relations} (that is, forgetting the reflection relations).
The finite and affine Weyl groups are Coxeter groups; we will make
precise the definition of the Artin groups in these cases.

\begin{Def} With the notation above define
\begin{enumerate}
\item[(i)] the finite Artin group $\A_{\W}$ as the group generated
by elements
$$T_1,\dots,T_n$$
satisfying the same braid relations as the reflections
$s_1,\dots,s_n$; \item[(ii)] the affine Artin group $\A_W$ as the
group generated by the elements
$$T_0,\dots,T_n$$
satisfying the same braid relations as the reflections
$s_0,\dots,s_n$.
\end{enumerate}
\end{Def}
From the definition it is clear that the finite Artin group can be
realized as a subgroup inside the affine Artin group. For further
use we also introduce the following lattices
$\cal{Q}_Y:=\{Y_\mu;\mu\in \Q\}$ and $\cal{Q}_X:=\{X_\b;\b\in
\Q\}$. Recall that the affine Weyl group also has a second
presentation, as a semidirect product. There is a corresponding
description of the affine Artin group  due to van der Lek
\cite{lek} and independently obtained by Lusztig \cite{lusztig}
and Bernstein (unpublished). To be more precise, Lusztig and
Bernstein give a proof for the corresponding description of the
{\sl extended Hecke algebra} (the proof also works for the
extended Artin group), which is easier and purely algebraic. H.
van der Lek's result is more difficult and the proof uses the
topological construction on the affine Artin group.
\begin{Prop}\label{lusztigpresentation}
The affine Artin group $\A_W$ is generated by the finite Artin
group and the lattice $\cal{Q}_Y$ such that the following
relations are satisfied for all $1\leq i\leq n$
\begin{enumerate}
\item[(i)] $T_iY_\mu=Y_\mu T_i \ \ \text{if \ }
(\mu,\a_i^\vee)=0$, \item[(ii)] $T_iY_\mu T_i=Y_{s_i(\mu)}   \ \
\text{if \ } (\mu,\a_i^\vee)=1$.
\end{enumerate}
\end{Prop}

\begin{Rem}\label{braid01}
In this description $Y_\mu=T_{\l_\mu}$ for $\mu$ any anti-dominant
element of the root lattice. For example
$Y_{-a_0^{-1}\th}=T_{s_\th}T_0$. In fact, the above Proposition
implies that the element $T_{s_\th}^{-1}Y_{-a_0^{-1}\th}(=T_0)$
satisfies the predicted braid relations with the generators $T_i$.
\end{Rem}
The special form of the Coxeter relations makes clear that the
affine Artin group admits an {\sl anti-involution} $\iota$ which
fixes all its generators. From Proposition
\ref{lusztigpresentation} and the above Remark it follows that the
element
$T_{s_\th}^{-1}\iota(Y_{-a_0^{-1}\th})=T_{s_\th}^{-1}T_0T_{s_\th}$,
and as a consequence also
$\iota(T_{s_\th}^{-1}T_0T_{s_\th})=T_{s_\th}T_0T_{s_\th}^{-1}$,
satisfies  with the generators $T_i$ the same braid relations as
$T_0$. This also follows from the topological description of the
affine Artin group.
\begin{Prop}\label{affineautos}
The cyclic group of infinite order $\Z$ acts on the affine Artin
group as automorphisms by conjugation with $T_{s_\th}$ on $T_0$
and by fixing the rest of the generators.
\end{Prop}

In fact van der Lek's description is even more precise; he
identifies a finite set of relations which should be imposed.
\begin{Prop}
The affine Artin group $\A_W$ is generated by the finite Artin
group and the lattice $\cal{Q}_Y$ such that the following
relations are satisfied for $1\leq i,j\leq n$
\begin{enumerate}
\item[(i)] For any pair of indices $(i,j)$ such that
  $2r_{ji}=-(\a_j,\a_i^\vee)$,
with $r_{ji}$ a non-negative integer  we have
\begin{equation}
T_iY_{\mu_j}=Y_{\mu_j} T_i
\end{equation}
where $\mu_j=\a_j+r_{ji}\a_i$; note that $(\mu_j,\a_i^\vee)=0$,
\item[(ii)] For any pair of indices $(i,j)$ such that \
$2r_{ji}-1=-(\a_j,\a_i^\vee)$, with $r_{ji}$ a non-negative integer we
have
\begin{equation}
T_iY_{\mu_j} T_i=Y_{s_i(\mu_j)}
\end{equation}
where $\mu_j=\a_j+r_{ji}\a_i$; note that $(\mu_j,\a_i^\vee)=1$.
\end{enumerate}
If we are in the $A_{2n}^{(2)}$ case the relations for the pair
$(i,n)$ are those obtained by considering $a_0^{-1}\a_n$ instead
of $\a_n$ in the above formulas.
\end{Prop}
Cherednik \cite{c1} (in the case of a reduced root system) and
Sahi \cite{sahi} (for a nonreduced root system) defined the
so-called double affine Hecke algebra. We extract here from their
definition only the definition of the double affine Artin group.
\begin{Def}\label{defcherednik}
The double affine Artin group $\A_{\tilde W}$ is generated by the
affine Artin group $\A_W$, the lattice $\cal{Q}_X$ and the element
$X_{a_0^{-1}\d}$ such the following relations are satisfied for
all \ $0\leq i\leq n$
\begin{enumerate}
\item[(i)] $X_{a_0^{-1}\d}$ is central\ , \item[(ii)]
$T_iX_\b=X_\b T_i \ \ \text{if \ } (\b,\a_i^\vee)=0$, \item[(iii)]
$T_iX_\b T_i=X_{s_i(\b)}   \ \ \text{if \ } (\b,\a_i^\vee)=-1.$
\end{enumerate}
\end{Def}
To avoid cumbersome notation we will make use of the convention
$Y_{a_0^{-1}\d}:=X_{-a_0^{-1}\d}$. The following fact exploits the
similarity of the above definition with presentation of the affine
Artin group given in Proposition \ref{lusztigpresentation}.
\begin{Rem}\label{braid03}
 The elements $X_{a_0^{-1}\th}T_{s_\th}^{-1}$ and $T_0$
satisfy the same braid relations with the generators $T_i$. In
fact, subgroup  of the double affine Artin group generated by
$T_i$, $i\neq 0$ and the lattice $\cal{Q}_X$ could be described as
follows. It is the group with generators $T_i$, $i\neq 0$ and
$X_{a_0^{-1}\th}$ and such that the elements $T_i$,
$X_{a_0^{-1}\th}T_{s_\th}^{-1}$ satisfy the same braid relations
as $s_i$, $i\neq 0$ and $s_0$. The relations (ii) and (iii) in the
above Proposition can be used to {define} generators for the
commutative lattice $\cal{Q}_X$.
\end{Rem}

The double affine Weyl group is not a Coxeter group, but a
generalized Coxeter group (in the sense of Saito and Takebayashi,
see \cite{saito2}) and we can define the associated Artin group in
the same way as for a Coxeter group (i.e. by keeping the
generalized braid relations and forgetting the reflection
relations). The equivalence between the two definitions has been
recently established in \cite{take}.

\begin{Def}
The elliptic Artin group is the factor group of the double affine
Artin group by the group generated by $X_{a_0^{-1}\d}$.
\end{Def}
As before, let us give a more refined description of the double
affine Artin group. It is based, as in the case of the affine
Artin group, on the topological description of the double affine
Artin group \cite{ion}.
\begin{Prop}\label{firstpresent}
The double affine Artin group $\A_{\tilde W}$ is generated by the
affine Artin group and the lattice $\cal{Q}_X$ such that the
following relations are satisfied for $0\leq i \leq n$ and $1\leq
j\leq n$
\begin{enumerate}
\item[(i)] For any pair of indices $(i,j)$ such that
  $2r_{ji}=-(\a_j,\a_i^\vee)$,
with $r_{ji}$ a non-negative integer  we have
\begin{equation}\label{eq1}
T_iX_{\mu_j}=X_{\mu_j} T_i
\end{equation}
where $\mu_j=\a_j+r_{ji}\a_i$; note that $(\mu_j,\a_i^\vee)=0$,
\item[(ii)] For any pair of indices $(i,j)$ such that \
$2r_{ji}+1=-(\a_j,\a_i^\vee)$, with $r_{ji}$ a non-negative integer we
have
\begin{equation}\label{eq2}
T_iX_{\mu_j} T_i=X_{s_i(\mu_j)}
\end{equation}
where $\mu_j=\a_j+r_{ji}\a_i$; note that $(\mu_j,\a_i^\vee)=-1$,
\item[(iii)] $X_{a_0^{-1}\d}$ is central.
\end{enumerate}
If we are in the $A_{2n}^{(2)}$ case the relations for the pair
$(i,n)$ are those obtained by considering $a_0^{-1}\a_n$ instead
of $\a_n$ in the above formulas.
\end{Prop}
To define the Hecke algebras, we introduce a field $\F$ (of
parameters) as follows: fix indeterminates $q$ and
$t_{0},\dots,t_n$ such that $ t_j=t_k$ if and only if $d_j=d_k;$
let $m$ be the lowest common denominator of the rational numbers
$\{(\a_j,\l_k)\ |\ 1\leq j,k\leq n \}$, and let $\F$ denote the
field of rational functions in $q^{\frac{1}{m}}$ and
$t_j^{\frac{1}{2}}$. For the Cartan matrix $A_{2n}^{(2)}$ we need
in fact more parameters: in this case two more indeterminates are
added to our field $t_{03}$ and $t_{02}$. To keep a uniform
notation we use $t_{01}$ to refer to $t_0$ in this case.
\begin{Def}
\noi \begin{enumerate}
\item[(i)] The finite Hecke algebra $\H_{\W}$ is the quotient of
the group $\F$-algebra of the finite Artin group by the relations
\begin{equation}
T_j-T_j^{-1}=t_j^{\frac{1}{2}} -t_j^{-\frac{1}{2}},
\end{equation}
for $1\leq j\leq n$. \item[(ii)] The affine Hecke algebra $\H_W$
is the quotient of the group $\F$-algebra of the affine Artin
group by the relations
\begin{equation}\label{q1}
T_j-T_j^{-1}=t_j^{\frac{1}{2}} -t_j^{-\frac{1}{2}},
\end{equation}
for all $0\leq j\leq n$.
\end{enumerate}
\end{Def}
Note here that by \cite[Remark 1.7]{ion} the relation (\ref{q1})
for $j=0$ needn't be imposed if the root system is reduced, since
it is a consequence of the other relations. However it is
absolutely necessary to impose it for $A_{2n}^{(2)}$. The
definition of the double affine Hecke algebra makes use of this
fact.
\begin{Def}\label{def3}
The double affine Hecke algebra $\H_{\tilde W}$ is the quotient of
the group $\F$-algebra of the double affine Artin group by the
relations
\begin{equation}\label{def1}
T_j-T_j^{-1}=t_j^{\frac{1}{2}} -t_j^{-\frac{1}{2}},\ \ \ \text{for
all } \ \ 1\leq j\leq n,
\end{equation}
and
\begin{equation}
X_\d=q^{-1}.
\end{equation}
If the root system is non-reduced the following relations are also
imposed
\begin{equation}
T_0-T_0^{-1}=t_{01}^{\frac{1}{2}} -t_{01}^{-\frac{1}{2}},
\end{equation}
\begin{equation}\label{t2}
X_{a_0^{-1}\th}T_{s_\th}^{-1}-T_{s_\th}X_{-a_0^{-1}\th}=t_{03}^{\frac{1}{2}}
-t_{03}^{-\frac{1}{2}},
\end{equation}
\begin{equation}\label{t3}
T_0^{-1}X_{\a_0}-X_{-\a_0}T_0=t_{02}^{\frac{1}{2}}
-t_{02}^{-\frac{1}{2}}.
\end{equation}
\end{Def}
The above definition of the double affine Hecke algebra for
$A_{2n}^{(2)}$ coincides indeed with the definition given by Sahi
\cite{sahi} after setting $t_{03}=u_n$, $t_{02}=u_0$ and making
use of Lemma 4.1 in \cite{sahi}.

%%%%%%%%%%%%%%%%%%%%%%%%%%%%%%%%%%%%%%%%%%%%%%%%%%%
%%%%%%%%%%%%%%%%%%%%%%%%%%%%%%%%%%%%%%%%%%%%%%%%%%%
%%%%%%%%%%%%%%%%%%%%%%%%%%%%%%%%%%%%%%%%%%%%%%%%%%%
\medskip
\subsection{The modular group and the braid group on three letters}
Let us consider the modular group $\rm SL(2,\Z)$ of two-by-two
matrices with integer entries and determinant one. By $B_3$ we
denote the braid group acting on three letters. With generators
and relations it has the following description: $B_3$ is the group
generated by $\mf{a, b}$ satisfying the relation
\begin{equation}
\mf{aba=bab}.
\end{equation}
Inside $\rm SL(2,\Z)$ consider the following elements
$$
u_{12}=\left(
\begin{array}{cc}
1 & 1\\
0 & 1
\end{array}\right)\ \ \ \ \text{and}\ \ \ \
u_{21}=\left(
\begin{array}{cc}
1  & 0\\
-1 & 1
\end{array}\right)\ .
$$
\begin{Lm}
There exists a surjective morphism of groups
$$\pi: B_3 \to {\rm SL(2,\Z)}$$ defined by $\pi(\mf a)=u_{12}$,
$\pi(\mf b)=u_{21}$. The kernel of this morphism is the subgroup
of the braid group spanned by $\mf{c^2}$, with $\mf{ c=(aba)^2}$
the generator of the center of $B_3$.
\end{Lm}
\begin{proof}
Indeed, a simple computation shows that $u_{12}$ and $u_{21}$
satisfy the braid relation and also $(u_{12}u_{21})^6={\rm I}$,
where $\rm I$ is the identity matrix. It is a straightforward
check that imposing the above two relations on the two generators
$u_{12}$ and $u_{21}$ constitutes an abstract description of the
modular group.
\end{proof}
%%%%%%%%%%%%%%%%%%%%%%%%%%%%%%%%%%%%%%%%%%%%%%%%%%%
%%%%%%%%%%%%%%%%%%%%%%%%%%%%%%%%%%%%%%%%%%%%%%%%%%%
%%%%%%%%%%%%%%%%%%%%%%%%%%%%%%%%%%%%%%%%%%%%%%%%%%%
\section{Automorphisms of triple
 groups}

\subsection{A few remarks }

We start by a few comments on the definition of the triple %%%%%%%%%%%%%%%%
% affine Artin
group. There are three copies of the affine Artin group
imbedded inside the triple %%%%%%%%%%%%%%%%%%%%%%%%%%%%%%%%%%%%%%%%%%%%%%%%
% affine Artin
 group and their
interaction needs to be elucidated.  Each of them contains a
commutative subgroup isomorphic with the lattice $\Q$, whose
elements  will be denoted by $\cal X^{0i}_\mu$.

We shall see how slightly more complicated relations also hold. We
approach first the case of a single affine bond. The following
Lemma will be very useful.
\begin{Lm}\label{trick}
Assume that two elements $p$ and $q$ satisfy the single lace
Coxeter relations with a third element $x$. Then the following are
equivalent
\begin{enumerate}
\item[(i)] The element $qp$ satisfies the double lace Coxeter
relation with $x$;

\item[(ii)] The element $p^{-1}qp$ satisfies the single lace
Coxeter relation with $x$.
\end{enumerate}
\noi Furthermore, if the above equivalent conditions also hold the
element $(qp)^{-1}p(qp)$ satisfies the single lace Coxeter
relation with $x$.
\end{Lm}
\begin{proof}
To prove the equivalence let us note that
\begin{eqnarray*}
p^{-1}q p\cdot x\cdot p^{-1}q p &=&  p^{-1}q x^{-1} p \cdot x\cdot q p \\
&=& p^{-1}q x^{-1} q^{-1}\cdot q p \cdot x\cdot q p\\
&=& p^{-1}x^{-1} q^{-1}\cdot x\cdot q p \cdot x\cdot q p
\end{eqnarray*}
and also
\begin{eqnarray*}
x\cdot p^{-1}q p\cdot x&=&xp^{-1}x^{-1}\cdot x\cdot qp\cdot
x\\
&=& p^{-1} x^{-1} p \cdot x\cdot qp\cdot x\\
&=& p^{-1} x^{-1}q^{-1} \cdot q p \cdot x\cdot qp\cdot x
\end{eqnarray*}
In consequence the left hand sides are equal if and only if the
right hand sides are equal. The equivalence is now clear.

Suppose now that the equivalent conditions also hold.  Since
$$
p\cdot p^{-1}qp = qp
$$
and $p$, $p^{-1}qp$ satisfy the single lace Coxeter relation with
$x$ and $qp$ satisfies the double lace Coxeter relation, by the
first part of the Lemma we obtain that
$$
(p^{-1}qp)^{-1} p(p^{-1}qp)= (qp)^{-1}p(qp)
$$
satisfies the single lace Coxeter relation with $x$.
\end{proof}
An immediate consequence of this principle is the following
\begin{Prop}\label{braidautom1}
If $l_0=1$ the elements
$$
\T_{02}^{-1}\T_{01}\T_{02} \ \ \ \text{and}\ \  \
\T_{03}^{-1}\T_{02}^{-1}\T_{03}\T_{02}\T_{03}
$$
satisfy the single lace Coxeter relation with $\D^k\T_\a\D^{-k}$,
for any integer $k$.
\end{Prop}
\begin{proof}
Denote  $\D^k\T_\a\D^{-k}$ by $x$. For the first element we apply
the second part of previous Lemma for $p=\T_{01}$ and
$q=\T_{01}\T_{02}\T_{01}^{-1}$ and for the second element we apply
the Lemma for $p=T_{03}$ and $q=T_{02}$. The verification of the
hypothesis is straightforward.
\end{proof}
Next we shift to the case of an affine double bond.
\begin{Lm}\label{ellbraid1}
If the affine node is connected to some other node with two laces
then the elements $\T_{0i}$ and $\T_{0j}$ commute with $\T_\a
\T_{0i}\T_{0j} \T_\a$ for any $1\leq i< j\leq 3$.
\end{Lm}
\begin{proof}
We prove our claim only for $\T_{02}$.

Now, $\T_{02}$ commutes with $\T_\a \T_{01}\T_\a^{-1}$ (by
(\ref{ellbraid})) and with $\T_\a \T_{02} \T_\a$ (by braid
relations), therefore commutes with their product.

Hence, we have proved that
\begin{eqnarray}\label{ee}
\T_\a \T_{01}\T_{02} \T_\a \T_{02} &=& \T_{02} \T_\a
\T_{01}\T_{02} \T_\a .
\end{eqnarray}

A similar argument shows that the rest of our statement is true.
\end{proof}
The following result uses these observations.
\begin{Prop}\label{braidautom2}
If $l_0=2$ the elements
$$
\T_{02}^{-1}\T_{01}\T_{02}\ \ \ \and\ \ \
\T_{03}^{-1}\T_{02}\T_{03}.
$$
satisfy with the $\T_i$ ($i\neq 0$) the same braid relations as
$\T_{01}$ does.
\end{Prop}
\begin{proof}
For the nodes not connected with the affine node the braid
relations are clear, following directly from the ones with the
$\T_{0i}$'s. The only thing to be checked is the double lace
Coxeter relation with $\T_\a$ (recall that we denote by $\a$ the
simple root corresponding to the node connected to the affine
node).
\begin{eqnarray*}
\T_\a \T_{02}^{-1}\T_{01}\T_{02}\T_\a
\T_{02}^{-1}\T_{01}\T_{02}&=&
\T_\a \T_{02}^{-1}\T_\a^{-1}\T_\a \T_{01}\T_{02}\T_\a \T_{02}^{-1}\T_{01}\T_{02}\\
&=& \T_\a \T_{02}^{-1}\T_\a^{-1}\T_{02}^{-1}\T_{01}\T_{02}\T_\a
\T_{01}\T_{02}\T_\a
\hspace{.3cm} \text{by (\ref{ee})}\\
&=& \T_{02}^{-1}\T_\a^{-1}\T_{02}^{-1}\T_\a \T_{01}\T_{02}\T_\a
\T_{01}\T_{02}\T_\a\\
&=& \T_{02}^{-1}\T_{01}\T_{02}\T_\a
\T_{02}^{-1}\T_{01}\T_{02}\T_\a \hspace{1.35cm} \text{by
(\ref{ee})}.
\end{eqnarray*}
The result for $\T_{03}^{-1}\T_{02}\T_{03}$ follows along the same
lines.
\end{proof}
Similar things happen with respect to relation (\ref{ellbraid}).
\begin{Prop}\label{ellbraidautom}
If $l_0=2$, the relation (\ref{ellbraid}) holds also for the
following pairs of elements $$(\T_{02},
\T_{02}^{-1}\T_{01}\T_{02}), \ (\T_{01},
\T_{03}^{-1}\T_{02}\T_{03}),\  (\T_{02}^{-1}\T_{01}\T_{02},
\T_{03})\ \text{and}\  (\T_{03}, \T_{03}^{-1}\T_{02}\T_{03})$$
\end{Prop}
\begin{proof}
As before we show this only for one of the above pairs. The
argument can be easily repeated to explain the result for any of
the remaining pairs.

The element $\T_{02}$ commutes with
$\T_\a^{-1}\T_{02}^{-1}\T_{\a}^{-1}$ (by the braid relations) and
with  $\T_\a\T_{01}\T_{02}\T_\a$ (by (\ref{ee})). Therefore it
commutes with their product. We have showed that
$$
\T_{02}\T_\a^{-1} \T_{02}^{-1}T_{01}\T_{02}\T_{\a}=\T_\a^{-1}
\T_{02}^{-1}T_{01}\T_{02}\T_{\a}\T_{02}
$$
which is precisely the desired relation.
\end{proof}
%%%%%%%%%%%%%%%%%%%%%%%%%%%%%%%%%%%%%%%%%%%%%%%%%%%%%%%%%%%%%%%%%%%%%%%%%%%%%
%%%%%%%%%%%%%%%%%%%%%%%%%%%%%%%%%%%%%%%%%%%%%%%%%%%%%%%%%%%%%%%%%%%%%%%%%%%%%
%%%%%%%%%%%%%%%%%%%%%%%%%%%%%%%%%%%%%%%%%%%%%%%%%%%%%%%%%%%%%%%%%%%%%%%%%%%%%
\subsection{An action of the braid group on three letters}
Next we will explain how the braid group on three letters act on
the triple %%%%%%%%%%%%%%%%%%%%%%%%%%%%%%%%%%%%%%%%%%%%%%%%%%%%%%%%%%%%%%%%%%%
% affine Artin
 group as automorphisms which fix all the
generators corresponding to the non-affine nodes. Before defining
this action let us construct an {\sl anti-automorphism} which
satisfies the same restrictions. Of course we only need to define
it only on $\T_{01}$, $\T_{03}$ and $\T_{02}$. Let $\mf e$ be the
anti-involution which fixes $\T_{02}$ and interchanges $\T_{01}$
and $\T_{03}$. This is of course possible since the defining
relations are invariant under such a change. It clear that this is
an involution, its square being the identity morphism.

Let us define next two special elements morphisms of the triple %%%%%%%%%%%%%%%%%%%%
group. The first one will be called $\underline{\mf a}$ and it
will fix $\T_{03}$ and it will send $\T_{01}$ to $\T_{02}$ and
$\T_{02}$ to $\T_{02}^{-1}\T_{01}\T_{02}$, and the second one will
be called $\underline{\mf b}$ and it will fix $\T_{01}$ and it
will send $\T_{02}$ to $\T_{03}$ and $\T_{03}$ to
$\T_{03}^{-1}\T_{02}\T_{03}$.
\begin{Thm}\label{tripleartinauto}
The above maps define indeed automorphisms of the triple %%%%%%%%%%%%%%%%%%%%%%%%%%%%
% affine Artin
group. Moreover, the map sending ${\mf a}$ and ${\mf b}$ to
$\underline{\mf a}$ and $\underline{\mf b}$, respectively, defines
a group morphism
$$
\Upsilon: B_3\to {Aut}(\A).
$$
\end{Thm}
\begin{proof}
Let us first argue that $\underline{\mf a}$ defines indeed a
endomorphism of the double affine Artin group. Indeed, by
Proposition \ref{braidautom1} and Proposition \ref{braidautom2}
the images of the generators $T_{0i}$ satisfy indeed the braid
relations.

In the case of a single affine bond we have to check that the
images of $\T_{03}^{-1}\T_{02}\T_{03}$ satisfy the required braid
relations. The image by ${\mf a}$ is $\D^{-1}\T_{01}\D$ and there
is nothing to prove and the image by ${\mf b}$ is
$\T_{03}^{-1}\T_{02}^{-1}\T_{03}\T_{02}\T_{03}$ and the check it
was done in Proposition \ref{braidautom1}.

In the case of the double affine bond we need to also check that
the images of the relations (\ref{ellbraid}) hold. This was proved
in Proposition \ref{ellbraidautom}. Therefore ${\mf a}$ and ${\mf
b}$ are indeed endomorphisms of the triple %%%%%%%%%%%%%%%%%%%%%%%%%%%%%%%%%
 group.

They are indeed isomorphisms since a simple check shows that
$$
\underline{\mf e}\underline{\mf a}\underline{\mf e}=\underline{\mf
b}^{-1}.
$$
As for the second part of our statement, this is again a
straightforward check.
\end{proof}
We will later prove that this morphism is injective.
%%%%%%%%%%%%%%%%%%%%%%%%%%%%%%%%%%%%%%%%%%%%%%%%%%%%%%%%%%%%%%%%%%%%%%%%%%%%%
%%%%%%%%%%%%%%%%%%%%%%%%%%%%%%%%%%%%%%%%%%%%%%%%%%%%%%%%%%%%%%%%%%%%%%%%%%%%%
%%%%%%%%%%%%%%%%%%%%%%%%%%%%%%%%%%%%%%%%%%%%%%%%%%%%%%%%%%%%%%%%%%%%%%%%%%%%%
\section{The triple %%%%%%%%%%%%%%%%%%%%%%%%%%%%%%%%%%%%%%%%%%%%%%%%%%%%%%%%%%%%%
% affine Artin
 group and the double affine Artin group}

\subsection{A refinement of Cherednik's presentation}
We start by analyzing the presentation of the double affine Artin
group in Proposition \ref{firstpresent}. We first focus on the
relations (\ref{eq1}) and (\ref{eq2}) for the pairs of type
$(0,j)$.
\begin{description}
\item[Type \ I] These are relations associated to $1\leq j\leq n$
such that $2r_{j0}=-(\a_j,\a_0^\vee)=(\a_j,\th)$, with $r_{j0}$ a
non-negative integer. For such a $j$ if we set $\mu_j=\a_j+r_{j0}\a_0$
the following relation holds
$$
T_0X_{\mu_j}=X_{\mu_j} T_0.
$$
Since $X_{a_0^{-1}\d}$ is central, we can certainly replace
$\mu_j$ in the above relation by $\a_j-r_{j0}a_0^{-1}\th$. The
only case in which the scalar product $(\a_j,\th)$ is even and
nonzero is for $A=B_n^{(2)}, A_{2n}^{(2)}$ when it could be equal
to $2$. Therefore, the relations of this type are
\begin{eqnarray}
T_0X_{\a_j} &=& X_{\a_j} T_0\ \ \ \text{if}\ \   (\a_j,\th)=0 \label{eq3} \\
T_0X_{\a_j-a_0^{-1}\th}&=&X_{\a_j-a_0^{-1}\th} T_0\ \ \ \text{if}\
\   (\a_j,\th)=2. \label{eq4}
\end{eqnarray}
Note the the relation (\ref{eq4}) is present only if $A=B_n^{(2)},
A_{2n}^{(2)}$. \item[Type II]  These are relations associated to
$1\leq j\leq n$ such that
$2r_{j0}+1=-(\a_j,\a_0^\vee)=(\a_j,\th)$, with $r_{j0}$ a non-negative
integer. For such a $j$ if we set $\mu_j=\a_j+r_{j0}\a_0$ the
following relation holds
$$
T_iX_{\mu_j} T_i=X_{s_i(\mu_j)}
$$
The only odd value the scalar product $(\a_j,\th)$ could take is
$1$, and this happens only if $A\neq B_n^{(2)}, A_{2n}^{(2)}$. In
the two excepted cases the scalar product takes only even values.
Therefore, the relations of this type are
\begin{equation}\label{eq5}
T_0X_{\a_j}T_0 = X_{\a_j+\a_0}\ \ \ \text{if}\ \ (\a_j,\th)=1 .
\end{equation}
Note the the relation (\ref{eq5}) is not present if $A=B_n^{(2)},
A_{2n}^{(2)}$.
\end{description}

\begin{Prop}\label{reduction2}
Assuming the notation above, we can reduce the number of relations
in the presentation of the double affine Artin group by keeping
from all relations of type I and II described above only the
following one:
\begin{enumerate}
\item[(i)] $ T_0X_{\a}T_0 = X_{\a+\a_0} \ \ \ \text{if}\ \ l_0=1$,
\item[(ii)] $ T_0X_{\a-a_0^{-1}\th}=X_{\a-a_0^{-1}\th} T_0\ \ \
\text{if}\ \   l_0=2$.
\end{enumerate}
\end{Prop}
\begin{proof}
Assume that $\b$ and $\gamma$ are simple roots
 whose nodes in the Dynkin diagram are connected, but none of them is
connected to the node of $\a_0$. We also assume that $\b$ is
shorter that $\gamma$ or they have the same length. This will
imply that $ (\b,\gamma^\vee)=-1 $ always and consequently
$$
s_\gamma(\b)=\b+\gamma.
$$
From equation (\ref{eq2}) we know that
$$
T_\gamma X_\b T_\gamma =X_{\a+\b}
$$
or equivalently $X_\gamma=T_\gamma X_\b T_\gamma X_{-\b}$. From
this expression it is clear that if we know that $T_0$ satisfies
the braid relations and it commutes with $X_\b$ it will follow
that $T_0$ commutes with $X_\gamma$.

Let $l_0=1$ and $\a_j$ a simple root for which $(\a_j,\th)=0$ (in
other words the nodes of $\a_0$ and $\a_j$ are not connected in
the Dynkin diagram). Using the above remark several times if
necessary we see that the relation (\ref{eq3}) is implied by the
knowledge of the braid relations and of the commutation of $T_0$
with $X_\b$, where $\b$ is any neighbor of $\a$. The commutation
of $T_0$ and $X_\b$ holds indeed since $\a$ is short (remember
that $l_0=1$) and therefore as explained above
$$
X_\b=T_\b X_\a T_\b X_{-\a}.
$$
Now,
\begin{eqnarray*}
T_0 X_\b T_0^{-1} &=& T_0 T_\b X_\a T_\b X_{-\a}T_0^{-1} \\
&=& T_\b T_0 X_\a T_0 T_\b T_0^{-1} X_{-\a}T_0^{-1}\ \  \ \text{by
the braid
relations for $T_0$ and $T_\b$}  \\
&=& T_\b X_{\a+\a_0}  T_\b X_{-\a-\a_0}\ \ \  \ \ \ \ \ \  \ \
\text{by the hypothesis} \\
&=& T_\b X_{\a}  T_\b X_{-\a} \ \ \ \  \ \ \ \ \ \ \ \ \ \ \ \ \ \
\ \ \text{by the commutation of
$T_\b$ and $X_{\a_0}$} \\
&=& X_\b .
\end{eqnarray*}
This computation shows that if we impose on $T_0$ the relation
stated in the hypothesis (besides the braid relations) then all
relations of type I automatically hold. In the case $A_n^{(1)}$
there are two type II relations: one which we have by hypothesis
and another one, associated to the second neighbor (let us call it
$\a^\prime$) of the affine simple root in the Dynkin diagram. A
straightforward computation, which exploits the fact that  $T_0$
commutes with $X_{\a+\a^\prime +\a_0}$, (fact which is a
consequence of the commuting relations proved above) will show
that the type II relation for $\a^\prime$ holds. For completeness,
let us explain the details:
\begin{eqnarray*}
T_0 X_{\a^\prime} T_0 &=& T_0 X_{-\a-\a_0} X_{\a+\a^\prime+\a_0} T_0 \\
&=& T_0 X_{-\a-\a_0} T_0 X_{\a+\a^\prime+\a_0}\ \  \
\text{by the commuting relations}  \\
&=& X_{-\a}  X_{\a+\a^\prime+\a_0}\ \ \  \ \ \ \ \ \  \ \ \ \ \
\text{by the hypothesis} \\
&=&  X_{\a^\prime+\a_0}.
\end{eqnarray*}
The proof of our result in the case $l_0=1$ is now completed. The
case $l_0=2$ is treated completely similarly.
\end{proof}
%%%%%%%%%%%%%%%%%%%%%%%%%%%%%%%%%%%%%%%%%%%%%%%%%%%
%%%%%%%%%%%%%%%%%%%%%%%%%%%%%%%%%%%%%%%%%%%%%%%%%%%
\subsection{Some relations}

We present  here  some relations which hold inside the double
affine Artin group. The relations will be useful later.
\begin{Prop}\label{braid02}
The elements $T_i$ ($i\neq 0$) and $T_0^{-1}X_{\a_0}$ satisfy
inside the double affine Artin group the same braid relations as
$T_i$ and $T_0$.
\end{Prop}
\begin{proof}
The claim is obvious if $(\th,\a_i)=0$, or equivalently if the
nodes in question are not connected by laces in the Dynkin
diagram. The other possible values for the scalar product are
$(\th,\a_i)=1$ (if there is a single lace connecting the $\a_0$
and $\a_i$; consequently $\a_i$ is short) and $(\th,\a_i^\vee)=2$
(if there are two laces connecting the $\a_0$ and $\a_i$;
consequently $\a_i$ is not short). We will consider them
separately. In what follows we use our convention to denote by
$\a$ the simple root whose node is connected to the affine node in
the Dynkin diagram.

First, if  $(\th,\a)=1$ by Definition \ref{defcherednik} we know
that
\begin{eqnarray}
X_{\a_0}T_\a&=&T_\a^{-1}X_{\a_0+\a}\ \ \ \text{and} \label{a} \\
X_{\a+\a_0}T_0^{-1}&=&T_0X_{\a}.  \label{aa}
\end{eqnarray}
Now,
\begin{eqnarray*}
T_0^{-1}X_{\a_0}T_\a T_0^{-1}X_{\a_0}&=&
T_0^{-1}T_\a^{-1}X_{\a_0+\a}T_0^{-1}X_{\a_0} \ \ \ \text{by (\ref{a})}\\
&=& T_0^{-1}T_\a^{-1}T_0 X_{\a}X_{\a_0} \hspace{1.15cm} \text{by (\ref{aa})}\\
&=& T_\a T_0^{-1}T_\a^{-1}X_{\a_0+\a} \hspace{1.2cm}
\text{by the braid relation}\\
&=& T_\a T_0^{-1}X_{\a_0}T_\a \hspace{1.85cm} \text{by (\ref{a})}.
\end{eqnarray*}
We proved the desired braid relation (recall that in this case the
relevant nodes in the Dynkin diagram are connected by a single
lace):
$$
T_0^{-1}X_{\a_0}T_\a T_0^{-1}X_{\a_0}=T_\a T_0^{-1}X_{\a_0}T_\a.
$$

Second, if  $(\th,\a)=2$ from Definition \ref{defcherednik} we
know that
\begin{eqnarray}
X_{\a_0}T_\a &=& T_\a^{-1}X_{\a_0+\a} \label{b} \\
X_{\a+\a_0}T_0^{-1}&=&T_0^{-1}X_{\a+\a_0}  \label{bb} \\
X_{\a+2\a_0}T_\a &=&T_\a X_{\a+2\a_0}.  \label{bbb}
\end{eqnarray}
In the same manner,
\begin{eqnarray*}
T_0^{-1}X_{\a_0}T_\a T_0^{-1}X_{\a_0}T_\a &=&
T_0^{-1}T_\a^{-1}X_{\a_0+\a}T_0^{-1}X_{\a_0}T_\a \hspace{.65cm}
\text{by (\ref{b})}\\
&=& T_0^{-1}T_\a^{-1}T_0^{-1}X_{\a_0+\a}X_{\a_0}T_\a
\hspace{.65cm}
\text{by (\ref{bb})}\\
&=& T_0^{-1}T_\a^{-1}T_0^{-1}T_\a X_{2\a_0+\a} \hspace{1.15cm}
\text{by (\ref{bbb})} \\
&=& T_\a T_0^{-1}T_\a^{-1}T_0^{-1}X_{2\a_0+\a} \hspace{1.15cm}
\text{by the braid relation}\\
&=& T_\a T_0^{-1}T_\a^{-1}X_{\a_0+\a}T_0^{-1}X_{\a_0}
\hspace{.65cm} \text{by (\ref{bb})}\\
&=& T_\a T_0^{-1}X_{\a_0}T_\a T_0^{-1}X_{\a_0}\hspace{1.3cm}
\text{by (\ref{b})}.
\end{eqnarray*}
We proved the desired braid relation (recall that in this case the
relevant nodes in the Dynkin diagram are connected by two laces):
$$
T_0^{-1}X_{\a_0}T_\a T_0^{-1}X_{\a_0}T_\a =T_\a
T_0^{-1}X_{\a_0}T_\a T_0^{-1}X_{\a_0}.
$$
The proof is completed.
\end{proof}
If $l_0=2$ another important relation holds.
\begin{Lm}
With the notation above, if the affine node is connected by a
double lace with the node corresponding to the simple root $\a$,
the elements $T_0$ and $T_\a^{-1}T_{0}^{-1}X_{\a_0}T_\a$ commute
inside the double affine Artin group.
\end{Lm}
\begin{proof} Indeed,
\begin{eqnarray*}
T_{0}T_\a^{-1}T_{0}^{-1}X_{\a_0}T_\a &=&
T_{0}T_\a^{-1}T_{0}^{-1}T_\a^{-1}X_{\a_0+\a}
\hspace{1cm}\text{by (\ref{b})}\\
&=& T_\a^{-1}T_{0}^{-1}T_\a^{-1}T_{0}X_{\a_0+\a}
\hspace{1cm} \text{by the braid relation} \\
&=& T_\a^{-1}T_{0}^{-1}T_\a^{-1}X_{\a_0+\a}T_{0}
\hspace{1cm}\text{by (\ref{bb})} \\
&=& T_\a^{-1}T_{0}^{-1}X_{\a_0}T_\a T_{0} \hspace{1.65cm}\text{by
(\ref{b})}
\end{eqnarray*}
\end{proof}
%%%%%%%%%%%%%%%%%%%%%%%%%%%%%%%%%%%%%%%%%%%%%%%%%%%%%%%%%%%%%%%%
%%%%%%%%%%%%%%%%%%%%%%%%%%%%%%%%%%%%%%%%%%%%%%%%%%%%%%%%%%%%%%%%
%%%%%%%%%%%%%%%%%%%%%%%%%%%%%%%%%%%%%%%%%%%%%%%%%%%%%%%%%%%%%%%

\subsection{A quotient of the triple
 group}

The following element will play a {\sl central} role. If, with the
usual notation, $s_\th=s_{j_1}\cdots s_{j_m}$ is a reduced
decomposition of the reflection $s_\th$ in terms of simple
reflections, we denote by $\T_{s_\th}$ the product $\T_{j_1}\cdots
\T_{j_m}$. Since we imposed the braid relations the definition of
$\T_\th$ will not depend on the reduced decomposition chosen.
Define
$$
\q:=\T_{01}\T_{02}\T_{03}\T_{s_\th}.
$$

Consider next the following quotient of the triple %%%%%%%%%%%%%%%%%%%%%%%%%%%%
% affine Artin
group.

\begin{Def}
Let $A$ be an irreducible affine Cartan matrix subject to our
restriction and $S(A)$ its Dynkin diagram. The group $\tilde \A$
is the quotient of the triple %%%%%%%%%%%%%%%%%%%%%%%%%%%%%%%%%%%%%%%%%%%%%%%%%%
% affine Artin
 group $\A$ by the
following relations
\begin{equation}\label{central}
\q:=\T_{01}\T_{02}\T_{03}\T_{s_\th} \ \ \ \text{is central}
\end{equation}
\end{Def}

From now on we will consider only the group $\tilde \A$, therefore
no confusion will arise if we denote the images of the triple %%%%%%%%%%%%%%%%%%
%affine Artin
group elements by the quotient map by the same
symbols.

Let us collect a few facts which immediately follow from the above
definition.
\begin{Lm}
Using the notation above, if  $l_0=1$, $\T_{02}$ could be
expressed in terms of the other generators.
\end{Lm}
\begin{proof}
Let us start with an immediate consequence of (\ref{central}).
Since from a straightforward check
$$
\T_{02}=\T_{\a}\T_{01}^{-1}\D\T_{03}^{-1}
\T_{\a}\T_{03}\D^{-1}\T_{01} \T_{\a}^{-1},
$$
using the centrality of $\q$ we get
$$
\T_{02}=\T_{\a}\T_{01}^{-1}\q^{-1}\D\T_{03}^{-1}
\T_{\a}\T_{03}\q\D^{-1}\T_{01} \T_{\a}^{-1},
$$
and using the relation (\ref{central}) we obtain the following
formula
\begin{equation}\label{magic1}
\T_{02}=\T_{01}^{-1}\T_{\a}^{-1}\T_{01}\T_{\a} \T_{s_\th}^{-1}
\T_{03}^{-1} \T_{\a}\T_{03} \T_{s_\th}\T_{01} \T_{\a}^{-1},
\end{equation}
which shows that indeed $\T_{02}$ is expressible in terms of the
other generators.
\end{proof}
\begin{Lm}\label{reducel=1}
In defining the group $\tilde \A$, if $l_0=1$, the relations in
Definition \ref{maindef} b) are superfluous.
\end{Lm}
\begin{proof}
Because $\D^{-k}\T_{\a}\D^{k}=\T_{s_\th}^k\T_\a \T_{s_\th}^{-k}$
the result follows from Proposition \ref{affineautos}. Also, by
the same principle, it is enough to prove that the simple lace
Coxeter relations between $\T_{03}^{-1}\T_{02}\T_{03}$ and $\T_\a$
hold. By the Lemma \ref{trick} this is equivalent with the double
lace Coxeter relations between $\T_{02}\T_{03}$ and $\T_\a$.

Let us consider the relation
$$
\T_{02}\T_{03}\T_\a\T_{02}\T_{03}\T_\a=\T_\a\T_{02}\T_{03}\T_\a\T_{02}\T_{03}.
$$
It essentially says that $\T_{02}\T_{03}$ commutes with
$\T_\a\T_{02}\T_{03}\T_\a$. But since
$$
\T_{02}\T_{03}=\q \T_{01}^{-1}\T_{s_\th}^{-1}
$$
and $\q$ is central this is equivalent to the fact that
$\T_{01}^{-1}\T_{s_\th}^{-1}$ and $\T_\a
\T_{01}^{-1}\T_{s_\th}^{-1}\T_\a$ commute. This is always true
since inside the affine Artin group generated by $\T_{01}$ and
$\T_i$, $i\neq 0$ since it is equivalent, using for example the
relations from Proposition \ref{lusztigpresentation} and Remark
\ref{braid01}, with the fact that $\cal X^{01}_{a_0^{-1}\th}$ and
$\cal X^{01}_{-\a+a_0^{-1}\th}$ commute. The proof is completed.
\end{proof}
\begin{Lm}\label{reducel=2}
In the above definition, if $l_0=2$, only one of the relations
(\ref{ellbraid}) should be imposed.
\end{Lm}
\begin{proof}
We will prove that if we impose only one relation in
(\ref{ellbraid}), say
\begin{equation}\label{c}
\T_{01}\T_\a^{-1}\T_{02}\T_\a=\T_\a^{-1}\T_{02}\T_\a\T_{01}
\end{equation}
the other two easily following from this one and  the fact that
$\q$ is central. We will illustrate this briefly.

The above relation says that $\T_{01}$ and
$\T_\a^{-1}\T_{02}^{-1}\T_\a$ commute. Since
 $\T_{01}$ also commutes with $\q$, $\T_\a^{-1}\T_{01}^{-1}\T_\a^{-1}$ (this
is just the braid relation) and $\T_\a\T_{s_\th}^{-1}\T_\a$ (this
is a relation which could be checked easily), it follows that
$\T_{01}$ commutes with their product, which is
$$
\T_\a^{-1}\T_{02}^{-1}\T_\a\T_\a^{-1}\T_{01}^{-1}\T_\a^{-1}\T_\a
\T_{s_\th}^{-1}\T_\a\q=\T_\a^{-1}\T_{03}\T_\a.
$$
We have just proved that
$$
\T_{01}\T_\a^{-1}\T_{03}\T_\a=\T_\a^{-1}\T_{03}\T_\a\T_{01}.
$$
The argument is the same if we choose to keep any other relation.
\end{proof}
%%%%%%%%%%%%%%%%%%%%%%%%%%%%%%%%%%%%%%%%%%%%%%%%%%%%%%%%%%%%%%%%%%%%%%%%%%%
%%%%%%%%%%%%%%%%%%%%%%%%%%%%%%%%%%%%%%%%%%%%%%%%%%%%%%%%%%%%%%%%%%%%%%%%%%%
%%%%%%%%%%%%%%%%%%%%%%%%%%%%%%%%%%%%%%%%%%%%%%%%%%%%%%%%%%%%%%%%%%%%%%%%%%%
\subsection{The main result} Our main goal is to prove that the
double affine Artin group $\A_{\tilde W}$ is isomorphic to the
group $\tilde \A$. Let us define the candidates for isomorphisms
between the two groups. Let the map
$$
\phi: {\tilde \A} \to \A_{\tilde W}
$$
be defined as the extension to a group morphism of the map
$$
\phi(\T_i)=T_i \ \ (\text{for }i\neq0), \ \ \ \phi(\T_{01})=T_0,\
\ \ \phi(\T_{03})=X_{a_0^{-1}\th}T_{s_\th}^{-1}, \ \ \
\phi(\T_{02})=T_0^{-1}X_{\a_0}.
$$
The next Proposition shows that this could indeed be done.
\begin{Prop}\label{lemma1}
The map $\phi: {\tilde \A} \to \A_{\tilde W}$ is well defined.
\end{Prop}
\begin{proof}
 From Remark \ref{braid03} and Proposition \ref{braid02} it follows that
the images of the generators of the group $\tilde \A$ satisfy the
required braid relations inside the double affine Artin group
$\A_{\tilde W}$. Also,
\begin{eqnarray*}
\phi(\q)&=&\phi(\T_{01}\T_{02}\T_{03}\T_{s_\th})\\
&=&T_0T_0^{-1}X_{\a_0}X_{a_0^{-1}\th}T_{s_\th}^{-1}T_{s_\th}\\
&=&X_{a_0^{-1}\d}
\end{eqnarray*}
which is central in the double affine Artin group. The only thing
which needs explanation is the fact that the image of relation
(\ref{ellbraid}) holds if $l_0=2$. In fact, as it follows from
Lemma \ref{reducel=2}, we need to do this only for one relation,
say
\begin{equation}
\T_{01}\T_\a^{-1}\T_{02}\T_\a=\T_\a^{-1}\T_{02}\T_\a\T_{01}.
\end{equation}

Indeed,
\begin{eqnarray*}
\phi(\T_{01}\T_\a^{-1}\T_{02}\T_\a)&=&
T_{0}T_\a^{-1}T_{0}^{-1}X_{\a_0}T_\a \\
&=& T_{0}T_\a^{-1}T_{0}^{-1}T_\a^{-1}X_{\a_0+\a}
\hspace{1cm}\text{by (\ref{b})}\\
&=& T_\a^{-1}T_{0}^{-1}T_\a^{-1}T_{0}X_{\a_0+\a}
\hspace{1cm} \text{by the braid relation} \\
&=& T_\a^{-1}T_{0}^{-1}T_\a^{-1}X_{\a_0+\a}T_{0}
\hspace{1cm}\text{by (\ref{bb})} \\
&=& T_\a^{-1}T_{0}^{-1}X_{\a_0}T_\a T_{0}
\hspace{1.65cm}\text{by (\ref{b})} \\
&=& \phi(\T_\a^{-1}\T_{02}\T_\a\T_{01}).
\end{eqnarray*}
The proof is completed.
\end{proof}
Let the map
$$
\psi: \A_{\tilde W} \to  {\tilde \A}
$$
be defined as the extension to a group morphism of the map
$$
\psi(T_i)=\T_i\ \ (\text{for }i\neq 0), \ \ \
\psi(T_{0})=\T_{01},\ \ \ \psi(X_{a_0^{-1}\th})=\T_{03}\T_{s_\th},
\ \ \ \psi(X_{a_0^{-1}\d})=\q.
$$

\begin{Prop}\label{lemma2}
The map $\psi: \A_{\tilde W} \to {\tilde \A}$ is well defined.
\end{Prop}
\begin{proof}
As noted in Remark \ref{braid03}, the elements for which we
defined $\psi$ are enough to generate the double affine Artin
group. This has the advantage of reducing the number of relations
for us to check to the images by $\psi$ of the  following
\begin{enumerate}
\item braid relations for $T_i$ ($i\neq 0$) and $T_0$; \item braid
relations for $T_i$ ($i\neq 0$) and
 $X_{a_0^{-1}\th}T_{s_\th}^{-1}$;
\item $X_{a_0^{-1}\d}$ is central; \item relations from
Proposition \ref{reduction2}.
\end{enumerate}
The only nontrivial check will be to prove that the images of the
relations stated in Proposition \ref{reduction2} hold inside the
group $\tilde \A$. Let us consider first the case when $l_0=1$. We
have to prove that
$$
\T_{01}\psi(X_\a)\T_{01}=\psi(X_{\a_0+a}).
$$
Since $T_\a X_{-\th}T_\a=X_{\a-\th}$ we know that
$$
X_\a= T_\a X_{-\th}T_\a X_\th .
$$
Therefore we want that
$$
\T_{01} \T_\a \T_{s_\th}^{-1} \T_{03}^{-1} \T_\a
\T_{03}\T_{s_\th}\T_{01}= \q \T_\a \T_{s_\th}^{-1} \T_{03}^{-1}
\T_{\a}
$$
or equivalently
$$
\T_{s_\th}\T_\a^{-1} \T_{01} \T_\a \T_{s_\th}^{-1} \T_{03}^{-1}
\T_\a \T_{03}\T_{s_\th}\T_{01} \T_{\a}^{-1}\T_{03}=\q.
$$
This immediately follows from equation (\ref{magic1}).

In the case $l_0=2$ we have to prove that
$$
\T_{01} \psi(X_{\a-a_0^{-1}\th})= \psi(X_{\a-a_0^{-1}\th})\T_{01}.
$$
As before,
$$
X_{\a-a_0^{-1}\th}= T_\a X_{-a_0^{-1}\th}T_\a  ,
$$
hence our statement is proved as follows
\begin{eqnarray*}
\T_{01} \psi(X_{\a-a_0^{-1}\th})&=&
\T_{01} \T_\a \T_{s_\th}^{-1}\T_{03}^{-1} \T_{\a}\\
&=&\q^{-1}\T_{01} \T_\a\T_{01}\T_{02}\T_{\a}
\hspace{1.55cm} \text{by (\ref{central})} \\
&=&\q^{-1}\T_{01} \T_\a\T_{01}\T_\a\T_\a^{-1}\T_{02}\T_{\a}
 \\
&=&\q^{-1}\T_\a\T_{01}\T_\a\T_{01}\T_\a^{-1}\T_{02}\T_{\a}
\hspace{.5cm} \text{by the braid relations}\\
&=&\q^{-1}\T_\a\T_{01}\T_\a\T_\a^{-1}\T_{02}\T_{\a}\T_{01}
\hspace{.5cm} \text{by (\ref{ellbraid})}\\
&=&\q^{-1}\T_\a\T_{01}\T_{02}\T_{\a}\T_{01}
\\
&=&\T_\a \T_{s_\th}^{-1}\T_{03}^{-1} \T_{\a}\T_{01}
\hspace{1.8cm} \text{by (\ref{central})}\\
&=&\psi(X_{\a-a_0^{-1}\th})\T_{01}.
\end{eqnarray*}
The proof of the Proposition is now complete.
\end{proof}
Our main result is the following.
\begin{Thm}\label{main}
The groups $\A_{\tilde W}$ and $\tilde \A$ are isomorphic.
\end{Thm}
\begin{proof}
The morphisms constructed in Proposition \ref{lemma1} and
Proposition \ref{lemma2} are inverses for each other, as we could
easily check this on generators.
\end{proof}
In other words, the definition of the group $\tilde \A$ could
serve as a definition for the double affine Artin group. We will
state this explicitly, by keeping only the non-redundant relations.
\begin{Thm}
Let $A$ be an irreducible affine Cartan matrix subject to our
restriction and $S(A)$ its Dynkin diagram. The double affine Artin
group $\A_{\tilde W}$ is given by generators and relations as
follows:

\underline{Generators}: one generator $T_i$ for each node, with
the exception of the affine node for which we have three
generators $T_{01}$, $T_{02}$ and $T_{03}$.

\underline{Relations}: a)
 Braid relations for each pair of generators
associated to any pair of distinct nodes (note that there are
three generators associated to the affine node).

\hspace{1.52cm} b) If there are double laces connecting the affine
node with the node $\a$ (i.e. $l_0=2$) the following relation also
holds
\begin{equation}
T_{01}T_\a^{-1}T_{03}T_\a=T_\a^{-1}T_{03}T_\a T_{01}
\end{equation}

\hspace{1.52cm} c) The element
\begin{equation}
X_{a_0^{-1}\d}:=T_{01}T_{02}T_{03}T_{s_\th} \ \ \ \text{is
central}
\end{equation}
The elliptic Artin group has the same description except for the
last relation which is replaced by
\begin{equation}
X_{a_0^{-1}\d}=1
\end{equation}
\end{Thm}
By Definition \ref{def3} we obtain a new description of the double
affine Hecke algebra. Note that the elements appearing in
(\ref{t2}) and (\ref{t3}) are $\T_{03}$ and $\T_{02}$,
respectively.

The same result at the level of Weyl groups is also of interest.
\begin{Thm}\label{wmain}
Let $A$ be an irreducible affine Cartan matrix subject to our
restriction and $S(A)$ its Dynkin diagram. The double affine Weyl
group ${\tilde W}$ is given by generators and relations as
follows:

\underline{Generators}: one generator $s_i$ for each node, with
the exception of the affine node for which we have three
generators $s_{01}$, $s_{02}$ and $s_{03}$.

\underline{Relations}: a)
 Braid relations for each pair of generators
associated to any pair of distinct nodes (note that there are
three generators associated to the affine node).

\hspace{1.52cm} b) If there are double laces connecting the affine
node with the node $\a$ (i.e. $l_0=2$) the following relation also
holds
\begin{equation}\label{wellbraid}
s_{01}s_\a s_{03}s_\a=s_\a s_{03}s_\a s_{01}.
\end{equation}

\hspace{1.52cm} c) All generators have order two.

\hspace{1.52cm} d) If, with the usual notation, the following
relation holds
\begin{equation}\label{wcentral}
\tau_{a_0^{-1}\d}:=s_{01}s_{02}s_{03}{s_\th} \ \ \ \text{is
central}
\end{equation}
The elliptic Weyl group has the same description except for the
last relation which is replaced by
\begin{equation}
\tau_{a_0^{-1}\d}=1
\end{equation}
\end{Thm}
%%%%%%%%%%%%%%%%%%%%%%%%%%%%%%%%%%%%%%%%%%%%%%%%%%%
%%%%%%%%%%%%%%%%%%%%%%%%%%%%%%%%%%%%%%%%%%%%%%%%%%%
%%%%%%%%%%%%%%%%%%%%%%%%%%%%%%%%%%%%%%%%%%%%%%%%%%%

\section{Automorphisms of double affine Artin groups}

\subsection{A reflection representation for the double affine
Weyl group}

The fact that the generators of the double affine Weyl group in
Theorem \ref{wmain} have order two suggests that it may have a
faithful representation in which the generators will act as
reflections. We construct here such a representation.

Assume we have fixed an irreducible affine Cartan matrix $A$ of
rank $n$ satisfying our assumptions. Let us define the following
$n+4$--dimensional real vector space
$$
V:=\Gc^* \oplus \Re\d_1\oplus \Re\d_2\oplus \Re\L_1\oplus \Re\L_2
$$
together with the nondegenerate bilinear form $(\cdot,\cdot)$
which extends the natural scalar product on $\Gc$
\begin{eqnarray*}
(\d_i,\L_j)&=& \d_{ij},\ \  i,j=1,2  \\
(\d_i,\a_k)&=& 0,\ \ \ \  i=1,2,\ \  k=1\cdots n\\
(\L_i,\a_k)&=& 0,\ \ \ \  i=1,2,\ \  k=1\cdots n .
\end{eqnarray*}
By $r$ we will denote the maximum number of laces in the Dynkin
diagram of $A$ and by $\RR_s$, $\RR_l$ the short, respectively
long, roots in $\RR$. The vector space $V$ contains the subset
$\tilde R$ defined as
$$\tilde R=(\RR_s+\Z\d_1+\Z\d_2)\cup (\RR_\ell+r\Z\d_1+r\Z\d_2),
\ \ \ \ \text{if $A\neq A^{(2)}_{2n}$},$$
$$\tilde R=(\RR_s+\Z\d_1+\Z\d_2)\cup (\RR_\ell+r\Z\d_1+r\Z\d_2)
\cup (\frac{1}{2}\RR_\ell+\Z\d_1+\Z\d_2+\frac{1}{2}(\d_1+\d_2)),
 \ \ \ \ \text{for $A^{(2)}_{2n}$}.$$
 For each $\tilde \a\in \tilde R$ and $v\in V$ define
$$
s_{\tilde \a}(v)= v-\frac{2(v,\tilde \a)}{(\tilde \a,\tilde
a)}\tilde \a .
$$
The $s_{\tilde \a}$ are reflections of $V$. We consider an action
of the double affine Weyl group $\tilde W$ on the vector space $V$
by letting the simple reflections $s_i$ acting as $s_{\a_i}$, and
$s_{01}$,  $s_{03}$ and  $s_{02}$ acting as
$s_{a_0^{-1}(\d_1-\th)}$, $s_{a_0^{-1}(\d_2-\th)}$ and
$s_{a_0^{-1}(\d_1+\d_2-\th)}$, respectively. Also, the element
$\tau_{a_0^{-1}\d}$ acts as
$$
\tau_{a_0^{-1}\d}(x):=x+(x,\d_2)a_0^{-1}\d_1-(x,\d_1)a_0^{-1}\d_2.
$$
\begin{Prop}
The above action of the generators of the double affine Weyl group
extend to a faithful representation
$$
\rho: \tilde{W}\to {\rm GL}(V).
$$
\end{Prop}
\begin{proof}
All the relations in the Theorem \ref{wmain} are easily verified.
Hence, we defined indeed a representation of the double affine
Weyl group.  To check that it is faithful we use Proposition
\ref{wdef} which has as a consequence the fact that any element of
$\tilde W$ can be put in the form $\w\l_{\mu}\tau_\b$ and the
following explicit actions of $\l_\mu$ and $\tau_\b$:
$$\l_\mu(x)=x-(x,\mu)\d_1
+(x,\d_1)(\mu-\frac{1}{2}|\mu|^2\d_1)$$
$$\tau_\b(x)=x-
(x,\b)\d_2 +(x,\d_2)(\b-\frac{1}{2}|\b|^2\d_2).$$ The verification
is straightforward.
\end{proof}
From now on we will not distinguish between the double affine Weyl
group and its image inside ${\rm GL}(V)$ via the above
representation.

Note the vector space $V_{(0,0)}:=\Gc^* \oplus \Re\d_1\oplus
\Re\d_2$ is invariant under the action of $\tilde W$ and that the
element $\tau_{a_0^{-1}\d}$ acts trivially on it. We immediately
obtain that the restriction to $V_{(0,0)}$ gives a faithful
representation of the elliptic Weyl group.

\begin{Lm}
All the reflections $s_{\tilde \a}$ for $\a\in \tilde R$ belong to
the double affine Weyl group.
\end{Lm}
\begin{proof}
This is a simple check, using the fact that
$$
ws_{\tilde \a}w^{-1}=s_{w(\tilde \a)}
$$
for any $w\in {\tilde W}$ and ${\tilde \a}\in {\tilde R}$, and the
fact that $s_{a_0^{-1}(\d_1-\th)}$ and $s_{a_0^{-1}(\d_2-\th)}$
belong to the double affine Weyl group.
\end{proof}
\noi We consider next the following subgroup of the full automorphisms
group of $\tilde W$
$$
{\mf G}({\tilde W}):=\{f\in Aut({\tilde W})\ |\ f(s_i)=s_i, \
i\neq 0, \ \ \text{and}\ \ f(\tau_{a_0^{-1}\d})=
\tau_{a_0^{-1}\d}\}.
$$

Since the generators $s_{a_0^{-1}(\d_1-\th)}$ and
$s_{a_0^{-1}(\d_2-\th)}$ could be replaced by
$s_{a_0^{-1}(\overline\d_1-\th)}$ and
$s_{a_0^{-1}(\overline\d_2-\th)}$ as long as $\overline \d_1$ and
$\overline \d_2$ generate the lattice $\Z\d_1\oplus\Z\d_2$, for
any matrix $u$ in ${\rm SL}(2,\Z)$ we can define such an
automorphism which sends $s_{01}$ and $s_{03}$ to
$s_{a_0^{-1}(u\cdot\d_1-\th)}$ and $s_{a_0^{-1}(u\cdot\d_2-\th)}$,
respectively, Here ${\rm SL}(2,\Z)$ acts as usually on the lattice
generated by $\d_1$ and $\d_2$. We will denote this automorphism
by $\underline u$. Let us observe that
$$
{\underline u}(s_{m\d_1+n\d_2-\th})=s_{u\cdot(m\d_1+n\d_2)-\th}.
$$ The following result is an immediate consequence of the above
remarks.
\begin{Prop}The map described above is an injective morphism
$$
\upsilon: {\rm SL}(2,\Z)\to {\mf G}({\tilde W}).
$$
\end{Prop}\label{wauto}
Let us note that the above morphism could be extended to ${\rm
GL(2,\Z)}$. The element $$ e:=\left(
\begin{array}{cc}
  0 & 1 \\
  1 & 0
\end{array}
\right)
$$
gives rise to an involution $\underline e$ of $\tilde W$ which
normalizes the group of automorphisms given by the modular group.
For example
$$eu_{12}e=u_{21}^{-1}.$$
The fact that such automorphisms exist follows automatically from
Theorem \ref{tripleartinauto} since the double affine Weyl group
is a quotient of the triple %%%%%%%%%%%%%%%%%%%%%%%%%%%%%%%%%%%%%%%%%%%%%%%%%
 group by relations which are fixed by
the action of $B_3$ and by the fact that $\mf{c^2}$ acts trivially
on this quotient. We preferred the explicit geometrical argument
since it shows that actually ${\rm SL}(2,\Z)$ acts faithfully on
the double affine Weyl group.

%%%%%%%%%%%%%%%%%%%%%%%%%%%%%%%%%%%%%%%%%%%%%%%%%%%%%%%%%%%%%%%%%%%%%%%%%%%%%
%%%%%%%%%%%%%%%%%%%%%%%%%%%%%%%%%%%%%%%%%%%%%%%%%%%%%%%%%%%%%%%%%%%%%%%%%%%%%
%%%%%%%%%%%%%%%%%%%%%%%%%%%%%%%%%%%%%%%%%%%%%%%%%%%%%%%%%%%%%%%%%%%%%%%%%%%%%

%%%%%%%%%%%%%%%%%%%%%%%%%%%%%%%%%%%%%%%%%%%%%%%%%%%%%%%%%%%%%%%%%%%%%%%%%%%%%
%%%%%%%%%%%%%%%%%%%%%%%%%%%%%%%%%%%%%%%%%%%%%%%%%%%%%%%%%%%%%%%%%%%%%%%%%%%%%
%%%%%%%%%%%%%%%%%%%%%%%%%%%%%%%%%%%%%%%%%%%%%%%%%%%%%%%%%%%%%%%%%%%%%%%%%%%%%

\subsection{Automorphisms of double affine Artin groups}
Returning to the double affine Artin group, denote by $T_{0i}$ the
images of the elements $\T_{0i}$ in $\A_{\tilde W}$. We will be
interested in automorphisms of the double affine Artin group which
descend to automorphism of the double affine Weyl group. Let us
call the group of such automorphisms $Aut(\A_{\tilde W};{\tilde
W})$. As for double affine Weyl groups, let us consider the
following group
$$
{\mf G}(\A_{\tilde W}):=\{f\in Aut(\A_{\tilde W};{\tilde W})\ |\
f(T_i)=T_i, \ i\neq 0, \ \ \text{and}\ \ f(X_{a_0^{-1}\d})=
X_{a_0^{-1}\d}\}.
$$

As we remarked in the case of the double affine Weyl group, the
action of the braid group on three letters described in Theorem
\ref{tripleartinauto} descends to the double affine Artin group,
since this is a factor of the triple %%%%%%%%%%%%%%%%%%%%%%%%%%%%%%%%%%%%%%%%%%%%
 group by a relation which is
invariant under $B_3$. It is obvious that these descents belong to
${\mf G}(\A_{\tilde W})$. Even more, the anti-involution $\mf e$,
which fixes $\T_{02}$ and interchanges $\T_{01}$ and $\T_{03}$,
also descends to the double affine Artin group. To keep our
notation simple and because no confusion will arise we will use
the same symbols to denote automorphisms of $\A_{\tilde W}$ which
descend from those of $\A$.

It easy to see that $\mf e$ interchanges $X_\mu$ with $Y_{-\mu}$
for any element of the root lattice. On any group we have a
canonical anti-involution which sends any element to its inverse.
If we consider its composition with $\mf e$ we obtain an
involution $\underline{\mf e}$ which could be described as
follows: it sends each of the generators $T_i$, $i\neq 0$ to its
inverse and interchanges $X_\mu$ with $Y_{\mu}$ for any element of
the root lattice. The descent of $\underline{\mf e}$ to the double
affine Weyl group is $\underline e$.

This involution -- or rather its descent to the double affine
Hecke algebra -- plays a central role in the theory of Macdonald
polynomials, where is responsible for the so-called difference
Fourier transform. To mention only one of its many implications we
note that the difference Fourier transform was the crucial
ingredient in Cherednik's proof \cite{c2} of the Macdonald
evaluation--duality conjecture. Let us state this result
explicitly.
\begin{Thm}\label{involution}
The map which sends each of the generators $T_i$, $i\neq 0$ to its
inverse and interchanges $X_\mu$ with $Y_{\mu}$ for any $\mu\in
\Q$ can be uniquely extended to an automorphism
$$
\underline{\mf e}: \A_{\tilde W}\to \A_{\tilde W}
$$
of the double affine Artin groups. The square of $\underline{\mf
e}$ is the identity isomorphism.
\end{Thm}
This a trivial consequence of our description of the double affine
Artin group, however it is quite hard to prove it starting from
the original description of Cherednik. See \cite[Chapter
3]{macbook}, \cite[Theorem 2.2]{ion} and \cite[Theorem 4.2]{sahi}
for a different proof. We restate Theorem \ref{tripleartinauto} in
the double affine Artin group context.

\begin{Thm}\label{doubleartinauto}
The map sending ${\mf a}$ and ${\mf b}$ to $\underline{\mf a}$ and
$\underline{\mf b}$, respectively, defines a group morphism
$$
\Upsilon: B_3\to {\mf G}(\A_{\tilde W}).
$$
\end{Thm}

By $w_\circ$ we denote the longest element of the the Weyl group
$\W$. The following result is well known.
\begin{Lm}
Inside the affine Artin group $\A_W$, $T_0$ commutes with
$T_{s_\th}^{-1}T_{w_\circ}$.
\end{Lm}
Next we  will describe the action of the center of the braid group
$B_3$ on the double affine Artin group.
\begin{Cor}\label{centeraction}
The generator $\mf c$ of the center of the braid group $B_3$ acts
(via $\Upsilon$) on the double affine Artin group  as the
conjugation by $T_{w_\circ}$.
\end{Cor}
\begin{proof}
By the above Theorem $\mf c$ acts as conjugation by $T_{s_\th}$ on
the elements $T_{0i}$. By the above Lemma this is the same as
conjugation by $T_{w_\circ}$.
\end{proof}

By Deligne \cite{deligne}, the center of the Artin group $\A_{\W}$
is generated by $T_{w_\circ}$ when $w_\circ=-1$ or by
$T_{w_\circ}^2$ if $w_\circ\neq -1$. The latter occurs only for
types $A_n$, $n\geq 2$, $D_{2n+1}$ and $E_6$.
\begin{Cor}\label{cheredthm}
The above action of $B_3$ will give a morphism from ${\rm
PSL}(2,\Z)$ (if $w_\circ=-1$) or ${\rm SL}(2,\Z)$ (if $w_\circ\neq
-1$) to the outer automorphism group of the double affine Artin
group.
\end{Cor}
The Corollary \ref{cheredthm}  recovers results of  Cherednik
\cite[Theorem 4.3]{c2}. The Theorem \ref{doubleartinauto} also
seems to follow from his work although no reference was made to
the braid group on three letters. The implications of the above
Corollary are extremely important. To give one example, it leads
to projective representations of the modular group expressed in
terms of special values of Macdonald polynomials. In turn these
give identities involving special values of Macdonald polynomials
at roots of unity.

Our approach allows us to say more about this action.
\begin{Thm}
The morphism $ \Upsilon: B_3\to {\mf G}(\A_{\tilde W}) $ is
injective.
\end{Thm}
\begin{proof} The following diagram is commutative
\begin{diagram}
B_3             & \rTo^{\Upsilon}  &  {\mf G}(\A_{\tilde W}) \\
\dTo^{\pi}      &                  &  \dTo                   \\
{\rm SL}(2,\Z)  & \rTo^{\upsilon}  &  {\mf G}({\tilde W})
\end{diagram}
Because $\upsilon$ is injective, any element in the kernel of
$\Upsilon$ must be contained in the kernel of $\pi$, which is
spanned by $\mf{c^2}$. From the above Corollary it is clear that
the image by $\Upsilon$ of such an element cannot act trivially on
${\mf G}(\A_{\tilde W})$. The Theorem is proved.
\end{proof}
This Theorem also implies that the original action of $B_3$ on the
triple %%%%%%%%%%%%%%%%%%%%%%%%%%%%%%%%%%%%%%%%%%%%%%%%%%%%%%%%%%%%%
 group is also faithful.
\begin{Thm}\label{triplefaith}
The morphism $ \Upsilon: B_3\to {Aut}(\A) $ is injective.
\end{Thm}

We close by noting that all the automorphisms of the double affine
Artin groups given by $B_3$ descend faithfully to the
corresponding double affine Hecke algebras, elliptic Artin groups
and elliptic Hecke algebras. For Hecke algebras there is no action
of $B_3$ on the parameters except for the non-reduced case where
it acts by permuting the $t_{0i}$.

%%%%%%%%%%%%%%%%%%%%%%%%%%%%%%%%%%%%%%%%%%%%%%%%%%%%%%%%%%%%%%%%%%%%%%%%%%%%%
%%%%%%%%%%%%%%%%%%%%%%%%%%%%%%%%%%%%%%%%%%%%%%%%%%%%%%%%%%%%%%%%%%%%%%%%%%%%%
%%%%%%%%%%%%%%%%%%%%%%%%%%%%%%%%%%%%%%%%%%%%%%%%%%%%%%%%%%%%%%%%%%%%%%%%%%%%%
\bibliographystyle{amsalpha}

\end{document}